\documentclass[a4paper,12pt]{amsart}
\usepackage{amssymb,verbatim}
\usepackage{ifthen}
\usepackage{graphicx}

\newtheorem{thm}[equation]{Theorem}
\newtheorem{cor}[equation]{Corollary}
\newtheorem{lem}[equation]{Lemma}

\newtheorem{rem}[equation]{Remark}

\newtheorem{conj}{Conjecture}
\theoremstyle{definition}
\newtheorem{example}[equation]{Example}
\newtheorem{prob}[equation]{Problem}

\newcounter {own}
\def\theown {\thesection       .\arabic{own}}

\newenvironment{pf}[1][]{%
 \vskip 3mm
 \noindent
 \ifthenelse{\equal{#1}{}}%
  {{\slshape Proof. }}%
  {{\slshape #1.} }%
 }%
{\qed\bigskip}

\newcounter{alphabet}
\newcounter{tmp}

\newcommand\Ref[1]{\setcounter{tmp}{1}\Alph{tmp}}

\newenvironment{Lem}[1][]{\refstepcounter{alphabet}%
\bigskip%
\noindent%
{\bf Lemma \Alph{alphabet}}%
{\bf .} \itshape}{\vskip 8pt}

\newcommand{\es}{{\mathcal S}}

\newcommand{\IC}{{\mathbb C}}
\newcommand{\ID}{{\mathbb D}}

\def\be{\begin{equation}}
\def\ee{\end{equation}}

\newcommand{\bee}{\begin{enumerate}}
\newcommand{\eee}{\end{enumerate}}

\newcommand{\blem}{\begin{lem}}
\newcommand{\elem}{\end{lem}}
\newcommand{\bthm}{\begin{thm}}
\newcommand{\ethm}{\end{thm}}
\newcommand{\bcor}{\begin{cor}}
\newcommand{\ecor}{\end{cor}}
\newcommand{\beg}{\begin{example}}
\newcommand{\eeg}{\end{example}}
\newcommand{\begs}{\begin{examples}}
\newcommand{\eegs}{\end{examples}}
\newcommand{\bdefe}{\begin{defin}}
\newcommand{\edefe}{\end{defin}}
\newcommand{\bprob}{\begin{prob}}
\newcommand{\eprob}{\end{prob}}
\newcommand{\bei}{\begin{itemize}}
\newcommand{\eei}{\end{itemize}}

\newcommand{\bcon}{\begin{conj}}
\newcommand{\econ}{\end{conj}}
\newcommand{\bcons}{\begin{conjs}}
\newcommand{\econs}{\end{conjs}}
\newcommand{\bprop}{\begin{propo}}
\newcommand{\eprop}{\end{propo}}
\newcommand{\br}{\begin{rem}}
\newcommand{\er}{\end{rem}}
\newcommand{\brs}{\begin{rems}}
\newcommand{\ers}{\end{rems}}
\newcommand{\bo}{\begin{obser}}
\newcommand{\eo}{\end{obser}}
\newcommand{\bos}{\begin{obsers}}
\newcommand{\eos}{\end{obsers}}
\newcommand{\bpf}{\begin{pf}}
\newcommand{\epf}{\end{pf}}
\newcommand{\ba}{\begin{array}}
\newcommand{\ea}{\end{array}}
\newcommand{\beq}{\begin{eqnarray}}
\newcommand{\beqq}{\begin{eqnarray*}}
\newcommand{\eeq}{\end{eqnarray}}
\newcommand{\eeqq}{\end{eqnarray*}}

\newcommand{\ds}{\displaystyle}

\def\cc{\setcounter{equation}{0}   
\setcounter{figure}{0}\setcounter{table}{0}}

\numberwithin{equation}{section}

\newcounter{minutes}
\divide\time by 60
\newcounter{hours}
\multiply\time by 60 \addtocounter{minutes}{-\time}

\begin{document}
\bibliographystyle{amsplain}
\title[Radius of Close-to-convexity of Harmonic Functions]{Radius of Close-to-convexity of Harmonic Functions}

\def\thefootnote{}
\footnotetext{ \texttt{\tiny File:~\jobname .tex,
          printed: \number\year-\number\month-\number\day,
          \thehours.\ifnum\theminutes<10{0}\fi\theminutes}
} \makeatletter\def\thefootnote{\@arabic\c@footnote}\makeatother

\author{David Kalaj}
\address{David Kalaj, University of Montenegro, faculty of natural sciences and mathematics,
Cetinjski put b.b. 81000, Podgorica, Montenegro}
\email{davidk@t-com.me}

\author{Saminathan  Ponnusamy}
\address{S. Ponnusamy, Department of Mathematics,
Indian Institute of Technology Madras, Chennai--600 036, India.}
\email{samy@iitm.ac.in}

\author{Matti Vuorinen}
\address{Matti Vuorinen Department of Mathematics, University of Turku, FIN-20014 Turku,
Finland.} \email{vuorinen@utu.fi}


\begin{abstract}
Let ${\mathcal H}$ denote the class of all normalized complex-valued harmonic
functions $f=h+\overline{g}$ in the unit disk ${\mathbb D}$, and let  $K=H+\overline{G}$
denote the harmonic Koebe function. Let $a_n,b_n, A_n, B_n$ denote the Maclaurin coefficients of
$h,g,H,G$, and
$${\mathcal F}=\{f=h+\overline{g}\in {\mathcal H}:\,|a_n|\leq A_n ~\mbox{ and }~|b_n|\leq B_n ~\mbox{ for }~ n\geq 1 \}.
$$
We show that the radius of univalence of the family ${\mathcal F}$
is $0.112903\ldots $. We also show that this number is also the radius of the starlikeness of ${\mathcal F}$.
Analogous results are proved for a subclass of the class of harmonic convex functions in ${\mathcal H}$. These results
are obtained as a consequence of a new coefficient inequality for certain class of harmonic close-to-convex functions.
Surprisingly, the new coefficient condition helps to improve Bloch-Landau constant for bounded harmonic mappings.
\end{abstract}
\thanks{ }

\maketitle \pagestyle{myheadings}
\markboth{D. Kalaj, S.Ponnusamy, and M. Vuorinen}{Harmonic Mappings }
\cc

\section{Introduction and Main Results}

Denote by ${\mathcal H}$ the class of all complex-valued harmonic
functions $f$ in the unit disk ${\mathbb D}=\{z \in {\mathbb C}:\,
|z|<1\}$ normalized by $f(0)=0=f_z(0)-1 $. Each $f$ can be
decomposed as $f=h+\overline{g}$, where $g$ and $h$ are analytic in
$\ID$ so that \cite{Clunie-Small-84,Duren:Harmonic}
\be\label{p1-11-eq1}
h(z)=z+\sum _{n=2}^{\infty}a_nz^n~\mbox{ and }~g(z)=\sum _{n=1}^{\infty}b_nz^n.
\ee
Let ${\mathcal S}_H$ denote the class of univalent and orientation-preserving
functions $f=h+\overline{g}$ in ${\mathcal H}$. Then the Jacobian of
$f$ is given by $J_f(z)=|h'(z)|^2 -|g'(z)|^2$. We note that if $f=h+\overline{g}
\in {\mathcal S}_H$ and  $g(z)\equiv 0$ in ${\mathbb D}$, then $f=h\in {\mathcal S}$,
where  ${\mathcal S}$ denotes the well-known class of normalized
univalent analytic functions in ${\mathbb D}$. A necessary and sufficient condition (see
\cite{Clunie-Small-84} or Lewy \cite{lewy-36}) for a harmonic
function $f$ to be locally univalent in $\ID$ is that $J_f(z) >0$ in
$\ID$. The function $\omega (z)=g'(z)/h'(z)$ denotes the complex
dilatation of $f$. Thus, for $f=h+\overline g\in {\mathcal S}_H$
with $g'(0)=b_1$ and $|b_1|<1$ (because $J_f(0)=1-|b_1|^2>0$), the
function
$$ F=\frac{f-\overline{b_1f}}{1-|b_1|^2}
$$
is also in ${\mathcal S}_H$. Thus, it is customary to restrict our
attention to the subclass
$${\mathcal S}_H^0=\{f\in{\mathcal S}_H:\,f_{\overline z}(0)=0\}.
$$
The family $\mathcal{S}_H^0$ is known to be compact. The
uniqueness result of the Riemann mapping theorem does not extend to
these classes of harmonic functions, \cite{Clunie-Small-84,Duren:Harmonic}.
Several authors have studied
the subclass of functions that map $\ID$ onto specific domains, eg.
starlike domains, convex and close-to-convex domains. Let ${\mathcal
S}_H^*$ (${\mathcal K}_H$, ${\mathcal C}_H$ resp.) consist of all
sense-preserving harmonic mappings $f=h+\overline{g}\in{\mathcal H}$
of $\ID$ onto starlike (convex, close-to-convex, resp.) domains.
Denote by ${\mathcal S}_H^{*0}$ (${\mathcal K}_H^{0}$, ${\mathcal
C}_H^{0}$ resp.) the class consists of those functions $f$ in
${\mathcal S}_H^*$ (${\mathcal K}_H$, ${\mathcal C}_H$ resp.) for
which $f_{\overline z}(0)=0$.

In \cite[Lemma~5.15]{Clunie-Small-84}, Clunie and Sheil-Small proved the
following result.

\begin{Lem}\label{uni-theo3}
If $h,g$ are analytic in $\mathbb D$ with $|h'(0)|>|g'(0)|$ and
$h+\epsilon g$ is close-to-convex for each $\epsilon$,
$|\epsilon|=1$, then $f=h+\overline{g}$ is close-to-convex in $\ID$.
\end{Lem}

This lemma has been used to obtain many important results. In the
case of ${\mathcal S}_H^0$, we have the harmonic Koebe function
$K=H+\overline{G}$ in ${\mathcal S}_H^0$, where
\be\label{p1-11-eq4}
H(z)=\frac{z-\frac{1}{2}z^2+\frac{1}{6}z^3}{(1-z)^3} ~\textrm{ and
}~ G(z)=\frac{\frac{1}{2}z^2+\frac{1}{6}z^3}{(1-z)^3}.
\ee

We see that the function $K$ has the dilatation $\omega (z)=z$ and
$K$ maps the unit disk $\ID$ onto the slit plane $\IC\backslash \{u+iv:\, u\leq -1/6,~v=0\}$. Moreover,
\[H(z)=z+\sum_{n=2}^\infty A_nz^n ~\textrm{ and }~ G(z)=\sum_{n=2}^\infty B_nz^n,
\]
where
\be\label{p1-11-eq2a} A_n=\frac{1}{6}(2n+1)(n+1)~\textrm{ and
}~ B_n=\frac{1}{6}(2n-1)(n-1), \quad n\geq 1.
\ee
A well-known coefficient conjecture of Clunie and Sheil-Small \cite{Clunie-Small-84}, is that
if $f=h+\overline{g} \in {\mathcal S}_H^0$ then the Taylor
coefficients of the series of $h$ and $g$ satisfy the inequality
\be\label{p1-11-eq2} |a_n|\le A_n  ~\mbox{ and }~ |b_n|\le B_n
\mbox{ for all $n\geq 1$}.
\ee
Although, the coefficients conjecture remains an open problem for the full class ${\mathcal S}_H^0$, the
same has been verified for certain subclasses, namely, the class ${\mathcal T}_H$
(see \cite[Section 6.6]{Duren:Harmonic}) of harmonic univalent typically real functions, the class of harmonic
convex functions in one direction, harmonic starlike functions in
${\mathcal S}_H^0$ (see  \cite[Section 6.7]{Duren:Harmonic}), and the class of harmonic close-to-convex functions
(see \cite{Wang-Liang}).


It is interesting to know to what extent do the conditions
(\ref{p1-11-eq2}) influence the univalency of the normalized
harmonic function $f(z)$ and of all of its partial sums, namely,
$f_n(z)$ and $f_{\overline{m}}(z)$, where
$$f_n(z) =h_n(z)+\overline{g_m(z)} ~\mbox{if $n\geq m$}; ~f_{\overline{m}}(z) =h_n(z)+\overline{g_m(z)}
~ \mbox{if $m\geq n$}.
$$
Here $h_n(z)$ and $g_m(z)$ represent the $n$-th section/partial sums
of $h$ and $g$ given by
$$h_n(z) = z+\sum_{k=2}^n a_kz^k ~\mbox{ and }~g_m(z) = \sum_{k=1}^m b_kz^k,
$$
respectively. According to our notation, the degree of the
polynomials  $f_n(z)$ and $f_{\overline{m}}(z)$ is $n$ if $n=m$.

\bthm\label{p1-11-Th1} Let $h$ and $g$ have the form {\rm
(\ref{p1-11-eq1})} and the coefficients of the series satisfy the
conditions {\rm (\ref{p1-11-eq2})}. Then $f=h+\overline{g}$
is close-to-convex $($univalent$)$, and starlike in the disk $|z|<r_S$, where
$$r_S=1+\frac{\sqrt{2}}{4}-\sqrt{\sqrt{2}+\frac{1}{8}}\approx 0.112903
$$
is the root of the quadratic equation
$$\sqrt{2}\, r^2-(1+2\sqrt{2})\,r +\sqrt{2} -1=0
$$
in the interval $(0,1)$. The result is sharp.
\ethm

The radii problems for various subclasses of univalent harmonic mappings are open \cite[Problem 3.3]{Bsh-Lyzzaik2011}
(see also \cite{Clunie-Small-84,Duren:Harmonic,Sh-Small90,RusSali89}). However, Theorem \ref{p1-11-Th1} quickly yields

\bcor
The radius of close-to-convexity and the radius of starlikeness for mappings in ${\mathcal S}_H^{*0}$
(resp. ${\mathcal C}_H^0$ and  ${\mathcal T}_H$) is at least $0.112903$.
\ecor

Under the hypotheses of Theorem \ref{p1-11-Th1}, all the partial sums of $f$ are
close-to-convex $($univalent$)$, and starlike in $|z|<r_S$. Similar comments apply
to the next two results.

Another well-known result due to Clunie and Sheil-Small
\cite{Clunie-Small-84} states that the coefficients of the series of
$h$ and $g$ of every convex function $f=h+\overline{g}\in {\mathcal
K}_H^0$ satisfy the inequalities
\be\label{p1-11-eq5}
|a_n|\le \frac{n+1}{2}  ~\mbox{ and }~ |b_n|\le\frac{n-1}{2}  \mbox{ for all $n\geq 1$}.
\ee
Equality occurs for the function $L=M+\overline{N}
\in {\mathcal K}_H^0$, where
\be\label{p1-11-eq6}
M(z)=\frac{1}{2} \left (\frac{z}{1-z}+ \frac{z}{(1-z)^2}\right )
~\textrm{ and }~N(z)=\frac{1}{2} \left (\frac{z}{1-z}- \frac{z}{(1-z)^2}\right ).
\ee
We observe that
$$L(z)={\rm Re\,}\left (\frac{z}{1-z}\right )+{\rm Im\,}\left (\frac{z}{(1-z)^2}\right )
=z+\sum_{n=2}^\infty \frac{n+1}{2} z^n -\overline{ \sum_{n=2}^\infty \frac{n-1}{2} z^n}.
$$
At this place it is worth recalling that the convexity  (resp. starlikness) property is not a hereditary
property in the harmonic case, unlike the analytic case. For instance, the convex function $L$ maps the
subdisk $|z|<r$ onto a convex domain for $r\leq \sqrt{2}-1$, but onto a non-convex domain for
$\sqrt{2}-1<r<1$.

\bthm\label{p1-11-Th2}
Let $h$ and $g$ have the form {\rm (\ref{p1-11-eq1})} and the coefficients of the series satisfy the
conditions {\rm (\ref{p1-11-eq5})}. Then $f=h+\overline{g}$
is close-to-convex $($univalent$)$, and starlike in the disk
$|z|<r_S$, where
$$r_S=1+\frac{\sqrt[3]{-18+\sqrt{330}}}{6^{2/3}}-
\frac{1}{\sqrt[3]{6(-18+\sqrt{330})}} \approx 0.164878
$$
is the real root of the cubic equation
$$2r^3-6r^2+7r-1=0
$$
in the interval $(0,1)$. The result is sharp.
\ethm

Theorem \ref{p1-11-Th2} easily gives the following corollary although Theorem \ref{p1-11-Th2} is much more stronger.

\bcor
The radius of close-to-convexity and the radius of starlikeness for convex mappings in ${\mathcal S}_H^{0}$
is at least $0.164878$.
\ecor


\bthm\label{p1-11-Th3}
Let $h$ and $g$ have the form {\rm (\ref{p1-11-eq1})} with $|b_1|=|g'(0)|<1$, and the coefficients of
the series satisfy the conditions
$$ |a_n| +|b_n|\leq c ~\mbox{ for all $n\geq 2$}.
$$
Then $f=h+\overline{g}$ is close-to-convex $($univalent$)$, and starlike in the disk $|z|<r_S$, where
$$r_S=1-\sqrt{\frac{c}{c+1-|b_1|}}.
$$
The result is sharp.
\ethm

Theorem \ref{p1-11-Th3} helps to improve the Bloch-Landau's theorem
for bounded harmonic functions. Consider the class ${\mathcal
B}_H^M$ of a harmonic mapping $f$ of the unit disk $\ID$ with
$f(0)=f_{\overline{z}}(0)=f_{z}(0)-1=0$, and $|f(z)|<M$ for $z\in
\ID$. There are two important constants one is relative to the
domain of the function while the other one, namely the Bloch
constant, is defined relative to the range. In \cite{HG}, authors
proved that if $f\in {\mathcal B}_H^M$ then $f$ is univalent in
$|z|<\rho_0$ and $f(|z|<\rho_0)$ contains a disk $|w|<R_0$, where
$$\rho _0 \approx \frac{1}{11.105M} ~\mbox{ and }~R_0=\frac{\rho _0}{2} \approx \frac{1}{22.21M}.
$$
Better estimates were given in \cite{DN,Armen-06,Liu1,Liu} and
later in \cite{CPW-BMMSS-11}, see Table \ref{table2} in which the functions $\phi $ and $\psi$
are explicitly given by
$$ \phi(x)= \frac{x}{\sqrt{2}(x^2+x-1)} ~\mbox{ and }~
\psi (x)= \frac{1}{\sqrt{2}}\left[1+ \left ( \frac{x^2-1}{x}\right
)\log \left ( \frac{x^2-1}{x^2+x-1}\right )\right].
$$
This result is the best known but not sharp.

The purpose the next theorem is to give a new proof of one of these
results. Indeed our method of proof is simple and improves the best
known result. In fact our distortion estimate for $f\in{\mathcal
B}_H^M$ provides the radius of close-to-convexity and the radius
starlikeness of ${\mathcal B}_H^M$.

\begin{thm}\label{p1-11-Th4}
Let $f\in{\mathcal B}_H^M$.
Then $f=h+\overline{g}$ is close-to-convex $($univalent$)$ in the disk $|z|<r_0$, where
$$r_S=1-\sqrt{\frac{4M}{4M+\pi}}
$$
and $f(\mathbb{D}_{r_{0}})$ contains a univalent disk of radius at least
$$R_S=r_S- \frac{4M}{\pi}\frac {r_S^2}{1-r_S}.
$$
\end{thm}

\begin{table}
\center{\Small
\begin{tabular}{|c|c|c||c|c|c|lccc|}
\hline
$M$ 
& $r_{}=\phi (8M/\pi)$ & $R_0=\psi (8M/\pi)$ &  $M$ 
& $ r_S $ & $R_S$\\
\hline
1& 0.22421   &  0.12629    & 1 &  0.251602  &  0.143904 \\
2&  0.11992  &  0.06367    & 2 &  0.152633  &  0.082622 \\
3&  0.08311  &  0.04328    & 3 &  0.109765  &   0.0580693 \\
\hline
\end{tabular}}
\bigskip
\caption{The left side columns refer to Theorem 4 in
\cite{CPW-BMMSS-11} and the right side columns refer to Theorem
\ref{p1-11-Th4}. \label{table2}}
\end{table}


\section{Useful Lemmas and their Proofs}

We need the following two lemmas to prove our main results.

\blem \label{p1-11-Lem1} Let $h$ and $g$ have the form {\rm
(\ref{p1-11-eq1})} with $|b_1|<1$, $f=h+\overline{g}$, and satisfy
the condition
\be\label{BP1Ieq3}
\sum _{n=2}^{\infty}n|a_n| + \sum
_{n=1}^{\infty}n|b_n| \leq 1.
\ee
Then $f\in {\mathcal C}_H^2$, where ${\mathcal C}_H^2=\{f \in {\mathcal S}_H:\, \mbox{$
|f_z(z)-1|<1-|f_{\overline{z}}(z)|$ in $\mathbb D$}\}.$ The bound in
{\rm (\ref{BP1Ieq3})} is sharp as the harmonic function
$$ f(z)= z+\sum _{n=2}^{\infty}\frac{\epsilon _n}{n}z^n +\sum _{n=1}^{\infty}\frac{\epsilon _n'}{n}\overline{z^n},
$$
for which $\sum _{n=2}^{\infty} |\epsilon _n| +\sum _{n=1}^{\infty}
|\epsilon _n'|=1$, shows. \elem \bpf In \cite{Hiroshi-Samy-2010}, it
was shown that ${\rm Re\,}f_z(z)>|f_{\overline{z}}(z)|$  whenever
(\ref{BP1Ieq3}) holds. The proof of this lemma follows from an easy
modification of the proof of the corresponding result from
\cite{Hiroshi-Samy-2010}. For the sake of completeness, we include
the detail. Note that the coefficient inequality implies that both
$h$ and $g$ are analytic in $\mathbb D$. Thus, $f=h+\overline{g}$ is
harmonic in $\ID$. Without loss of generality, we may assume that
$f$ is not affine. Then, as $f_z=h'$ and
$f_{\overline{z}}=\overline{g'}$, it follows from the hypotheses
that
\begin{align*}
|h'(z) -1| & \leq   \sum _{n=2}^{\infty}n |a_n| \, |z|^{n-1}\\
& \leq  \sum _{n=2}^{\infty}n|a_n| \leq 1- \sum _{n=1}^{\infty}n|b_n|\\
& \leq 1- |g'(z)|
\end{align*}
implying that $f\in {\mathcal C}_H^2$ (since strict inequality
occurs either at the second  or fourth inequality). In particular,
${\rm Re}\, h'(z)>|g'(z)|$ in $\ID$ and hence, $f$ is locally
univalent in $\ID$. \epf

For example, the functions
$$f_n(z)=z+\frac{n+1}{2n^2}z^{n} +\frac{n-1}{2n^2}\overline{z^{n}} ~ \mbox{ for $n\geq 2$}
$$
satisfy the condition (\ref{BP1Ieq3}) and hence, belong to the class
${\mathcal C}_H^2$. In the following lemma, we show that functions
in ${\mathcal C}_H^2$ are indeed close-to-convex in $\ID$.

\blem\label{p1-11-Lem3} Let $h$ and $g$ have the form {\rm
(\ref{p1-11-eq1})} with $|b_1|<1$, $f=h+\overline{g}$. Suppose $f\in
{\mathcal C}_H^2$. Then, we have the following
\begin{enumerate}
\item[\textbf{(a)}] $f$ is close-to-convex in $\mathbb D$.
\item[\textbf{(b)}] $\big |\,|a_n|-|b_n|\,\big |\leq 1/n$ for $n\geq 2$ whenever $b_1=0$. The equality occurs,
for example, for the function
$$f(z)=z+ \frac{e^{i\theta}}{n}z^{n} ~\mbox{ or }~f(z)=z+ \frac{e^{i\theta}}{n}\overline{z^{n}}
~ \mbox{ for $n\geq 2$ and $\theta$ real}.
$$
\item[\textbf{(c)}]
$\ds \sum_{n=2}^{\infty}n^{2}(|a_{n}|^{2}+|b_{n}|^{2})\leq
1-|b_1|^2. $
\end{enumerate}
\elem \bpf First we prove part \textbf{(a)}. Let
$f=h+\overline{g}\in {\mathcal C}_H^2$ and $F=h+\epsilon g$, where
$|\epsilon |=1$. Then,
$$|F'(z)-1|\leq |h'(z)-1| +|g'(z)|<1
$$
showing that $F$ is analytic and close-to-convex in $\mathbb D$.
According to Lemma \Ref{uni-theo3}, it follows that the harmonic
function $f$ is also close-to-convex (and univalent) in $\mathbb D$.

Next, set $\omega (z)=F'(z)-1$. Then, as $b_1=g'(0)=0$, we have
$\omega(0)=0$ and $|\omega (z)|<1$ for $z\in \ID$. It is well-known
property that the coefficients of such an analytic function $\omega$
satisfy the inequality $|\omega ^{(n)}(0)|\leq n!$ for each $n\geq
1$. This gives the estimate
$$\left |na_n+\epsilon nb_n\right |\leq 1 ~\mbox{ for each $n\geq 2$}.
$$
As $|\epsilon |=1$, triangle inequality gives the proof for part
\textbf{(b)}.

For the proof of part \textbf{(c)}, we observe that
$$|F'(z)-1|= \left |\sum_{n=2}^{\infty}na_{n}z^{n-1} +\epsilon \sum_{n=1}^{\infty}nb_{n}z^{n-1}\right |<1,\quad z\in\ID.
$$
Therefore, with $z=re^{i\theta}$ for $r\in (0,1)$ and $0\leq  \theta
\leq 2\pi$, the last inequality gives
$$\sum_{n=2}^{\infty}n^{2}(|a_{n}|^{2}+|b_{n}|^{2})r^{2(n-1)} +|b_1|^2=
\frac{1}{2\pi}\int_{0}^{2\pi} |F'(re^{i\theta})-1|^{2}\,d\theta \leq
1.
$$
Letting $r\rightarrow 1^{-}$, we obtain the  inequality
$$\sum_{n=2}^{\infty}n^{2}(|a_{n}|^{2}+|b_{n}|^{2})\leq 1-|b_1|^2
$$
and the proof is complete. \epf

In \cite{Hiroshi-Samy-2010}, under the hypotheses of Lemma
\ref{p1-11-Lem1}, it was actually shown that $f\in {\mathcal
C}_H^1$, where
$${\mathcal C}_H^1=\{f \in {\mathcal S}_H:\, \mbox{${\rm Re\,} f_z(z)>|f_{\overline{z}}(z)|$ in $\mathbb D$}\}.
$$
Clearly, Lemma \ref{p1-11-Lem1} improves this result because of the
strict inclusion  ${\mathcal C}_H^2\subsetneqq {\mathcal C}_H^1$.
Later, in \cite{Bhara-samy-pre11}, it was also shown that if
$b_1=g'(0)=0$, then the coefficient condition (\ref{BP1Ieq3})
ensures that $f\in\es_H^{*0}$ (see also \cite{Silverman-98}). In view of
Lemma \ref{p1-11-Lem1}, the result of \cite{Bhara-samy-pre11,Silverman-98}  may be stated in
an improved form.

\blem\label{p1-11-Lem2}
Let $h$ and $g$ have the form {\rm (\ref{p1-11-eq1})} with  $b_1=g'(0)=0$, $f=h+\overline{g}$, and
satisfy the condition
\be\label{BP1Ieq6}
\ds \sum _{n=2}^{\infty}n|a_n| + \sum _{n=2}^{\infty}n|b_n| \leq 1.
\ee
Then $f\in  {\mathcal C}_H^2\cap \es_H^{*0}.$
\elem

The following generalization of Lemma \ref{p1-11-Lem1} is easy to
obtain and so we omit its details.

\bcor\label{p1-11-cor1} Let $h$ and $g$ have the form {\rm
(\ref{p1-11-eq1})} with $|b_1|<1-\beta$ for some $\beta \in [0,1)$,
and $f=h+\overline{g}$. Then we have the following:
\begin{enumerate}
\item[\textbf{(a)}] If the coefficients of $h$ and $g$ satisfy the condition
\be\label{p1-11-eq7}
\sum _{n=2}^{\infty}n|a_n| + \sum _{n=1}^{\infty}n|b_n| \leq 1-\beta ,
\ee
then $f\in {\mathcal C}_H^2(\beta)$, where
$${\mathcal C}_H^2(\beta)=\{f \in {\mathcal S}_H:\, \mbox{$ |f_z(z)-1|<1-\beta-|f_{\overline{z}}(z)|$ in $\mathbb D$}\}.
$$
In particular, $f$ is close-to-convex in $\ID$. The bound here is
sharp as the harmonic function
$$ f(z)= z+\sum _{n=2}^{\infty}\frac{\epsilon _n}{n}z^n
+\sum _{n=1}^{\infty}\frac{ \epsilon _n'}{n}\overline{z^n},
$$
for which $\sum _{n=2}^{\infty} |\epsilon _n| +\sum _{n=1}^{\infty}
|\epsilon _n'|=1-\beta$, shows.
\item[\textbf{(b)}] If $f\in {\mathcal C}_H^2(\beta)$, then one has
$$\big |\,|a_n|-|b_n|\,\big |\leq (1-\beta)/n ~\mbox{ for $n\geq 2$ whenever $b_1=0$}.
$$
The equality occurs, for example, for the function
$$f(z)=z+ (1-\beta)\frac{e^{i\theta}}{n}z^{n} ~\mbox{ or }~f(z)=z+ (1-\beta)\frac{e^{i\theta}}{n}\overline{z^{n}}
~ \mbox{ for $n\geq 2$ and $\theta$ real}.
$$
We also have
$$\ds \sum_{n=2}^{\infty}n^{2}(|a_{n}|^{2}+|b_{n}|^{2})\leq (1-\beta)^2-|b_1|^2.
$$
\end{enumerate}
\ecor

It is a matter of routine checking to see that the coefficient condition
(\ref{p1-11-eq7}) is necessary for $f=h+\overline{g}$ to belong to
${\mathcal C}_H^2(\beta)$ 
whenever the Taylor coefficients $a_n\leq 0$ for all $n\geq 2$, and
$b_n\leq 0$ for all $n\geq 1$.

\section{Proofs of Main Theorems}

\noindent \textbf{Proof of Theorem \ref{p1-11-Th1}.} Let $h$ and $g$
have the form (\ref{p1-11-eq1}) satisfying the coefficient
conditions (\ref{p1-11-eq2}). First we observe that $b_1=g'(0)=0$.
The conditions (\ref{p1-11-eq2}) implies that the series
(\ref{p1-11-eq1}) are convergent in the unit disk $|z|<1$, and
hence, the sum $h$ and $g$ are analytic in $\ID$. Thus,
$f=h+\overline{g}$ is harmonic in $\ID$. Let $0<r<1$, we let
$$f_r(z):=r^{-1}f(rz)=r^{-1}h(rz) +r^{-1}\overline{g(rz)}$$
so that $f_r(z)=h_r(z) +\overline{g_r(z)}$ and
$$f_r(z)=z+\sum _{n=2}^{\infty}a_n r^{n-1} z^{n} +\overline{\sum _{n=2}^{\infty}b_n r^{n-1} z^{n}}, ~z\in \ID.
$$
By hypotheses,   $|a_n|\leq A_n$ and $|b_n|\leq B_n$ for $n\geq 2$, where $A_n$ and $B_n$ are given by (\ref{p1-11-eq2a}).
Using these coefficient estimates, we obtain
\beqq
S&=&\sum _{n=2}^{\infty}n|a_n| r^{n-1} +\sum _{n=2}^{\infty}n |b_n| r^{n-1}\\
& \leq&  \sum _{n=2}^{\infty}nA_{n} r^{n-1} +\sum _{n=2}^{\infty}nB_{n} r^{n-1}.
\eeqq
We show that $f_r\in {\mathcal C}_H^2\cap \es_H^{*0}.$
According to Lemma \ref{p1-11-Lem2}, it suffices to show that $S\leq 1$. By the last inequality, $S\leq 1$
if $r$ satisfies the inequality
$$ \sum _{n=2}^{\infty}nA_{n}r^{n-1} \leq 1-\sum _{n=2}^{\infty}nB_{n}r^{n-1},
$$
or equivalently (as $A_{n}+B_n=(2n^2+1)/3$),
\be\label{p1-11-eq3}
2\sum _{n=2}^{\infty}n^3r^{n-1} + \sum _{n=2}^{\infty}nr^{n-1} \leq 3 .
\ee
As
$$\frac{r}{(1-r)^2}=\sum _{n=1}^{\infty}nr^{n} ~\mbox{ and }~\frac{r(1+r)}{(1-r)^3}=\sum _{n=1}^{\infty}n^2r^{n},
$$
it follows that
$$\frac{(1-r)(1+2r)+3r(1+r)}{(1-r)^4}=\sum _{n=1}^{\infty}n^3r^{n-1}
$$
and  (\ref{p1-11-eq3}) reduces to the inequality,
$$\frac{2(r^2+4r+1)}{(1-r)^4}+\frac{1}{(1-r)^2} \leq 6, ~\mbox{ i.e. }~ 2(1-r)^4-(1+r)^2\geq 0.
$$
This gives
$$\sqrt{2}(1-r)^2-(1+r)=\sqrt{2}\, r^2-(1+2\sqrt{2})\,r +\sqrt{2} -1\geq 0.
$$

Thus, from Lemma \ref{p1-11-Lem2}, $f_r$ is close-to-convex (univalent) in $\ID$ and
starlike in $\ID$ for all $0<r\leq r_S,$ where $r_S$ is the root of
the quadratic equation
$$\sqrt{2}\, r^2-(1+2\sqrt{2})\,r +\sqrt{2} -1=0
$$
in the interval $(0,1)$.  
In particular, $f$ is close-to-convex (univalent) and starlike in $|z|<r_S$.

\begin{figure}
\begin{center}
\includegraphics[width=10cm]{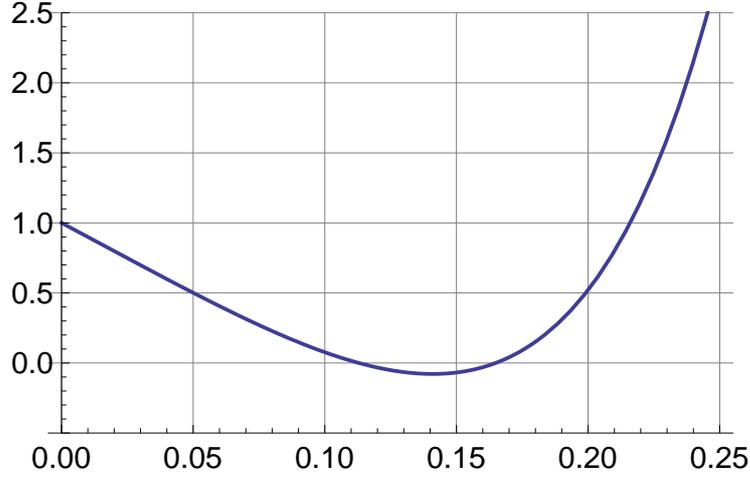}
\end{center}
\caption{The graph of the Jacobian $J_{F_0}(r)$ for $r\in (0,0.25)\,.$\label{fig1}}
\end{figure}

Next, to prove the sharpness part of the statement of the theorem,
we consider the function
$$F_0(z)= H_0(z)+\overline{G_0(z)}
$$
with
$$H_0(z)=2z-H(z) ~\mbox{ and }~ G_0(z)=-\overline{G(z)}.
$$
Here $H$ and $G$ are defined by (\ref{p1-11-eq4}). We note that
$$F_0(z)= z-\sum_{n=2}^\infty A_nz^n -\overline{\sum_{n=2}^\infty B_nz^n}.
$$
As $F_0$ has real coefficients we obtain.
\begin{align*}
J_{F_0}(r)
&=\left (H_0'(r) +G_0'(r)\right ) \left (H_0'(r) -G_0'(r)\right ) \\
&= \left(1- \sum _{n=2}^{\infty}nA_{n}r^{n-1} -\sum
_{n=2}^{\infty}nB_{n}r^{n-1}\right )
\left(1-\sum _{n=2}^{\infty}n(A_{n}-B_n) r^{n-1} \right )\\
&= \left(1- \sum _{n=2}^{\infty}\frac{n(2n^2+1)}{3} r^{n-1}  \right
) \left(1-\sum _{n=2}^{\infty}n^2 r^{n-1} \right ) \\&=\left(1 -
\frac{-4 r^2 + 3 r^3 - r^4}{(-1 + r)^3 r}\right) \left(1 + \frac{-6
r^2 + 5 r^3 -
    4 r^4 + r^5}{(-1 + r)^4 r}\right)\\&=\frac{(-1 + 7 r - 6 r^2 + 2 r^3) (1 - 10 r + 11 r^2 - 8 r^3 +
   2 r^4)}{(-1 + r)^7}.
\end{align*}
Thus $J_{F_0}(r)=0$, $0< r<1$ if and only if
$$r=r_S=\frac{1}{4} \left(4 + \sqrt{2} - \sqrt{2 + 16 \sqrt{2}}\right)\approx0.112903
$$
or
$$r=r'_S=1 + \left (-18 + \sqrt{330}\right )^{1/3}6^{-2/3} - \left (6 \big (-18 +
\sqrt{330}\big )\right )^{-1/3}\approx 0.164878.$$
Moreover for $r_S<r<r'_S$ we have $J_{F_0}(r)<0$. The graph of the function $J_{F_0}(r)$ for $r\in (0,0.25)$ is shown
in Figure \ref{fig1}.

 This observation together with Lewy's theorem
gives that (as the Jacobian changes sign), the function $F_0(z)$ is
not univalent in $|z|<r$ if $r> r_{\mathcal{S}}$ and thus,
$r_{\mathcal{S}}$ cannot be replaced by a larger number. \hfill
$\Box$

\vspace{8pt} \noindent \textbf{Proof of Theorem \ref{p1-11-Th2}.}
Following the notation and the method of the proof of Theorem
\ref{p1-11-Th1},  it suffices to show that $f_r\in {\mathcal C}_H^2\cap \es_H^{*0}$. According to
Lemma \ref{p1-11-Lem2},  $f_r\in {\mathcal C}_H^2\cap \es_H^{*0}$ whenever $S\leq 1$, where
$$S=\sum _{n=2}^{\infty}n|a_n| r^{n-1} +\sum _{n=2}^{\infty}n |b_n| r^{n-1}
$$
when $a_n$ and $b_n$ satisfy the coefficient inequalities given by \eqref{p1-11-eq5}. Finally,  using
\eqref{p1-11-eq5}, we see that $S\leq 1$ if $r$ satisfies the inequality
$$ \sum _{n=2}^{\infty}\frac{n(n+1)}{2}r^{n-1} \leq 1-\sum _{n=2}^{\infty}\frac{n(n-1)}{2}r^{n-1}.
$$
The last inequality is easily seen to be equivalent to
$$\frac{1}{2}\left [\frac{1}{(1-r)^2}+ \frac{1+r}{(1-r)^3} -1 \right ]\leq
1+ \frac{1}{2}\left [\frac{1}{(1-r)^2}- \frac{1+r}{(1-r)^3} -1
\right ]
$$
which upon simplification reduces to
$$2(1-r)^3-1-r=-(2r^3-6r^2+7r-1)\geq 0.
$$
The first part of the conclusion easily follows as in the proof of
Theorem \ref{p1-11-Th1}.

\begin{figure}
\begin{center}
\includegraphics[width=10cm]{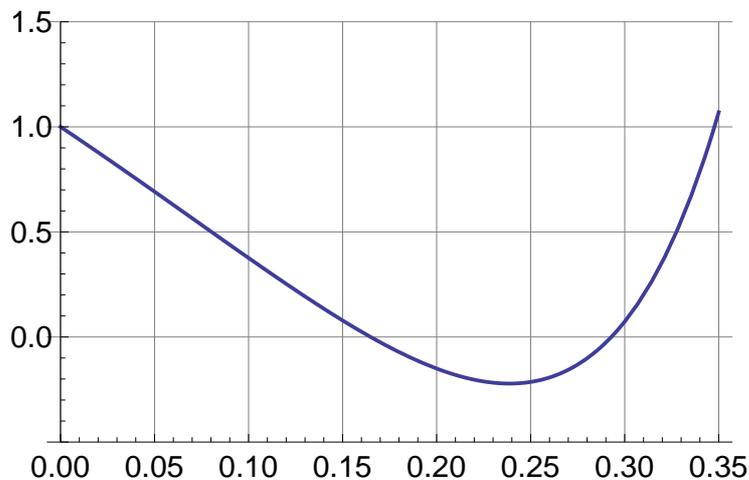}
\end{center}
\caption{The graph of the Jacobian $J_{L_0}(r)$ for $r\in (0, 0.35)\,.$\label{fig2}}
\end{figure}

The sharpness part of the statement of Theorem \ref{p1-11-Th2}
follows if we consider the function
$$L_0(z)=2z-M(z)-\overline{N(z)},
$$
where $M$ and $N$ are defined by (\ref{p1-11-eq6}).  We note that
$$L_0(z)= z-\sum_{n=2}^\infty \frac{n+1}{2}z^n +\overline{\sum_{n=2}^\infty \frac{n-1}{2}z^n}.
$$
Again, as $L_0$ has real coefficients, we can easily obtain that for $r\in (0,1)$
\begin{align*}
J_{L_0}(r) &=\left (2-M'(r)\right )^2- \left (N'(r)\right )^2 \\
&=\left (2-M'(r) +N'(r)\right ) \left ((2-M'(r) -N'(r)\right ) \\
&=   \left(2- \frac{1+r}{(1-r)^3}\right ) \left(2- \frac{1}{(1-r)^2} \right )\\
&=   \frac{2}{(1-r)^5}\left(2(1-r)^3- (1+r)\right ) \left(r-1- \frac{\sqrt2}{2} \right )\left(r-1+ \frac{\sqrt2}{2} \right ).
\end{align*}
We see that  $J_{L_0}(r_S)=0$, $0< r<1$ if and only if
$$r=r_S \approx 0.16487$$ or
$$r=r'_S= \frac{2-\sqrt2}{2} \approx0.292893 .
$$
Moreover for $r_S<r<r'_S$ we have $J_{L_0}(r)<0$. The graph of the function $J_{L_0}(r)$ for $r\in (0,0.35)$ is shown
in Figure \ref{fig2}.

Thus, according to Lewy's theorem, $L_0(z)$ is
not univalent in $|z|<r$ if $r> r_{\mathcal{S}}$ and this
observation shows that $r_{\mathcal{S}}$ cannot be replaced by a
larger number. \hfill $\Box$

\vspace{8pt}

\noindent \textbf{Proof of Theorem \ref{p1-11-Th3}.}
This time we apply Lemma \ref{p1-11-Lem1} and show that $f_r$ defined by
$f_r(z):=r^{-1}f(rz)=r^{-1}h(rz) +r^{-1}\overline{g(rz)}$ belongs to ${\mathcal C}_H^2$.

As in the proof of previous two theorems, it suffices to show the corresponding coefficient
inequality  \eqref{BP1Ieq3}, namely,
$$S= \sum _{n=2}^{\infty}n(|a_n| +|b_n|)r^{n-1}+|b_1|\leq 1.
$$
By the hypothesis, $ |a_n| +|b_n|\leq c$ for all $n\geq 2$ and so,
the last inequality $S\leq 1$ clearly holds if $r$ satisfies the inequality
$$ c\left (\frac{1}{(1-r)^2} -1 \right )\leq 1-|b_1|, ~\mbox{ i.e. }~r\leq r_S=1-\sqrt{\frac{c}{c+1-|b_1|}}.
$$
Thus, by Lemma \ref{p1-11-Lem1},
$$|h_r'(z)-1|<1-|g_r'(z)|
$$
holds for all $z\in\ID$ whenever $r\leq r_S$.  Thus,  $f\in {\mathcal C}_H^2$.

The function $f_0(z)=h_0(z)+\overline{g_0(z)}$, where
$$h_0(z)=z-\frac{c}{2} \left (\frac{z^2}{1-z}\right )  \textrm{ and }~ g_0(z)=-|b_1|z-\frac{c}{2} \left (\frac{z^2}{1-z}\right ),
$$
shows that the result is sharp. Indeed, it is easy to compute that
$$J_{f_0}(r)=|h_0'(r)|^2 -|g_0'(r)|^2=(1+|b_1|)\left (1+c-|b_1|-\frac{c}{(1-r)^2} \right )
$$
which shows that  $J_{f_0}(r_S)=0$ and $J_{f_0}(r)<0$ for  $r>r_S$. The proof of the theorem is complete.
\hfill $\Box$

\vspace{8pt}

\noindent \textbf{Proof of Theorem \ref{p1-11-Th4}.} Let
$f=h+\overline{g}$ be a harmonic mapping defined on the unit disk $\ID$ with
$f(0)=f_{\overline{z}}(0)=f_{z}(0)-1=0$, and $|f(z)|<M$ for $z\in
\ID$, where $h$ and $g$ have the form {\rm (\ref{p1-11-eq1})} with
$b_1=0$. According to \cite[Lemma 1]{CPW-JMAA-11} (see also
\cite{CPW-BMMSS-11}), we obtain the sharp estimates
\be\label{eq1.12c}
|a_{n}|+|b_{n}|\leq\frac{4M}{\pi} ~\mbox{  for any $n\geq 1$}.
\ee
As $b_1=0$ and $a_1=1$, it follows that $M\geq \pi /4 \approx 0.785398$. By Theorem \ref{p1-11-Th3}
with $c=4M/\pi$, we conclude that $f$ is close-to-convex and starlike (because $b_1=0$) for
$|z|<1-\sqrt{c/(c+1)}=r_S$.

In particular, $f$ is univalent for $|z|<r_S$ and furthermore, we have for $|z|=r_S$,
\beqq
|f(z)| &= &\left |z+\sum _{n=2}^{\infty}\big  (a_n z^n + \overline{b_nz^n}\big  )\right |\\
&\geq  &|z|-\left |\sum _{n=2}^{\infty}\big (a_n z^n + \overline{b_nz^n}\big ) \right |\\
&\geq  & r_S- \sum _{n=2}^{\infty}(|a_n| +|b_n|)r_S^n\\
&\geq   & r_S- \frac{4M}{\pi}\sum _{n=2}^{\infty}r_S^n\\
&=  & r_S- \frac{4M}{\pi}\frac {r_S^2}{1-r_S}=R_S
\eeqq
and the proof is complete.  \hfill $\Box$

\end{document}